\documentclass[letterpaper,10pt,reqno]{amsart}
\usepackage{indentfirst} 
\usepackage{amssymb}
\usepackage{mathtools} 
\usepackage{mathabx} 
\usepackage{amsthm}
\usepackage{thmtools}
\usepackage{enumitem} 
\usepackage[backref=page,colorlinks=true]{hyperref} 
\usepackage[usenames,dvipsnames]{xcolor} 
\usepackage{ifthen}
\usepackage{tikz}
\usetikzlibrary{decorations.pathreplacing}
\usepackage{tikz-cd}
  
\renewcommand*{\backref}[1]{}
\renewcommand*{\backrefalt}[4]{\tiny
  \ifcase #1 (\textbf{NOT CITED.})%
  \or    (Cited on page~#2.)%
  \else   (Cited on pages~#2.)%
  \fi}

\makeatletter

\def\MRbibitem{\@ifnextchar[\my@lbibitem\my@bibitem}

\def\mybiblabel#1#2{\@biblabel{{\hyperref{http://www.ams.org/mathscinet-getitem?mr=#1}{}{}{#2}}}}

\def\myhyperanchor#1{\Hy@raisedlink{\hyper@anchorstart{cite.#1}\hyper@anchorend}}

\def\my@lbibitem[#1]#2#3#4\par{%
  \item[\mybiblabel{#2}{#1}\myhyperanchor{#3}\hfill]#4%
  \@ifundefined{ifbackrefparscan}{}{\BR@backref{#3}}%
  \if@filesw{\let\protect\noexpand\immediate
    \write\@auxout{\string\bibcite{#3}{#1}}}\fi\ignorespaces%
}

\def\my@bibitem#1#2#3\par{%
  \refstepcounter\@listctr
  \item[\mybiblabel{#1}{\the\value\@listctr}\myhyperanchor{#2}\hfill]#3%
  \@ifundefined{ifbackrefparscan}{}{\BR@backref{#2}}%
  \if@filesw\immediate\write\@auxout
    {\string\bibcite{#2}{\the\value\@listctr}}\fi\ignorespaces%
}

\makeatother

\DeclareFontFamily{U} {MnSymbolA}{}
\DeclareFontShape{U}{MnSymbolA}{m}{n}{
   <-6> MnSymbolA5
   <6-7> MnSymbolA6
   <7-8> MnSymbolA7
   <8-9> MnSymbolA8
   <9-10> MnSymbolA9
   <10-12> MnSymbolA10
   <12-> MnSymbolA12}{}
\DeclareFontShape{U}{MnSymbolA}{b}{n}{
   <-6> MnSymbolA-Bold5
   <6-7> MnSymbolA-Bold6
   <7-8> MnSymbolA-Bold7
   <8-9> MnSymbolA-Bold8
   <9-10> MnSymbolA-Bold9
   <10-12> MnSymbolA-Bold10
   <12-> MnSymbolA-Bold12}{}
\DeclareSymbolFont{MnSyA} {U} {MnSymbolA}{m}{n}
 \DeclareFontFamily{U} {MnSymbolC}{}
\DeclareFontShape{U}{MnSymbolC}{m}{n}{
  <-6> MnSymbolC5
  <6-7> MnSymbolC6
  <7-8> MnSymbolC7
  <8-9> MnSymbolC8
  <9-10> MnSymbolC9
  <10-12> MnSymbolC10
  <12-> MnSymbolC12}{}
\DeclareFontShape{U}{MnSymbolC}{b}{n}{
  <-6> MnSymbolC-Bold5
  <6-7> MnSymbolC-Bold6
  <7-8> MnSymbolC-Bold7
  <8-9> MnSymbolC-Bold8
  <9-10> MnSymbolC-Bold9
  <10-12> MnSymbolC-Bold10
  <12-> MnSymbolC-Bold12}{}
\DeclareSymbolFont{MnSyC} {U} {MnSymbolC}{m}{n}

\DeclareMathSymbol{\top}{\mathord}{MnSyA}{219} 
\DeclareMathSymbol{\plus}{\mathord}{MnSyC}{20} 

\declaretheorem[numberwithin=section]{theorem}
\declaretheorem[sibling=theorem]{lemma}

\declaretheorem[sibling=theorem,style=definition]{definition}
\declaretheorem[sibling=theorem,style=remark]{example}

\declaretheorem[sibling=theorem,style=remark]{remark}
\declaretheorem[name=Acknowledgements, style=remark, numbered=no]{ack}

\hypersetup{bookmarksdepth = 3} 
\numberwithin{equation}{section}     

\setlist[enumerate,1]{label={\upshape(\alph*)},ref=\alph*}
\setlist[enumerate,2]{label={\upshape(\arabic*)},ref=\arabic*}

\newcommand{\R}{\mathbb{R}}

\newcommand{\N}{\mathbb{N}}

\DeclareMathOperator{\per}{per}

\newcommand{\Rnon}{\mathbb{R}_{\plus}}    
\newcommand{\Rpos}{\mathbb{R}_{\plus\plus}} 
\newcommand{\Mat}[2][]{\ifthenelse{\equal{#1}{}}{\R^{{#2}\times{#2}}}{\R^{{#1}\times{#2}}}}
\newcommand{\Man}[2][]{\ifthenelse{\equal{#1}{}}{\Rnon^{{#2}\times{#2}}}{\Rnon^{{#1}\times{#2}}}}
\newcommand{\Map}[2][]{\ifthenelse{\equal{#1}{}}{\Rpos^{{#2}\times{#2}}}{\Rpos^{{#1}\times{#2}}}}

\renewcommand{\epsilon}{\varepsilon}
\renewcommand{\phi}{\varphi}

\begin{document}

\title[Perfect matching in inhomogeneous random bipartite graphs]{Perfect Matchings in Inhomogeneous Random Bipartite Graphs in  Random Environment}
\date{\today}

\subjclass[2010]{Primary 05C80; Secondary 05C70, 05C63.}

\begin{thanks}
{J.B.\ was partially supported by Proyecto Fondecyt 1140202.
G.I.\ was partially supported by Proyecto Fondecyt 1150058.
M.P.\ was partially supported by Proyecto Fondecyt 1140988.}
\end{thanks}

\author[J.~Bochi]{Jairo Bochi}
\author[G.~Iommi]{Godofredo Iommi}
\author[M.~Ponce]{Mario Ponce}

\begin{abstract}
In this note we  study inhomogeneous  random bipartite graphs in random environment. These graphs can be thought of  as an extension of the classical  Erd\"os-R\'enyi random graphs in a random environment. We show that the expected number of perfect matchings obeys a  precise quenched asymptotic and that it can be approximated using an iterative process that converges exponentially fast.
\end{abstract}
\maketitle
\section{Introduction}
In their seminal paper \cite{erdos_Renyi}, Erd\"os and R\'enyi  studied the following random graphs that now bear their names.
Consider  a bipartite graph with set of vertices given by  $W=\left\{w_1, \dots, w_n   \right\}$ and $M=\left\{m_1, \dots, m_n   \right\}$.
Let $p \in [0,1]$ and consider the independent random variables $X_{(ij)}$ with law
\begin{eqnarray*}
X_{(ij)}(x)=
\begin{cases}
1 & \text{with probability } p ; \\
0 & \text{with probability } 1-p.
\end{cases}
\end{eqnarray*}
Denote by $G_n(x)$ the bipartite graph with vertices $W$ and $M$ and edges $E(x)$, where the edge $(w_i , m_j)$ belongs to $E(x)$ if and only if $X_{(ij)}(x)=1$.
Let $\text{pm}(G_n(x))$ be the number of \emph{perfect matchings} of the graph $G_n(x)$ (see Sec.~\ref{sec:4} for precise definitions). Erd\"os and R\'enyi \cite[p.460]{erdos_Renyi} observed that the mean of the number of perfect matchings is given by
\begin{equation} \label{er_formula}
\mathbb{E}(\text{pm}(G_n(x)))= n! p^n.
\end{equation}
This number has been also studied by Bollob\'as and McKay \cite[Theorem 1]{bolloMck}  in the context of $k-$regular random graphs and by O'Neil \cite[Theorem 1]{o} for random graphs having a fixed (large enough) proportion of edges. We refer to the text by Bollob\'as \cite{BOLLOBAS} for further details on the subject of random graphs.

This paper is devoted to study certain sequences of inhomogeneous  random bipartite graphs  $G_{n, \omega}$ in a random environment  $\omega \in \Omega$ (definitions are given in Sec.~\ref{sec:rbgre}). Inhomogeneous random graphs have been intensively studied over the last years  (see  \cite{BOLLOBAS2}).  Our main result  (see Theorem~\ref{thm:main} for precise statement) is that there exists a constant $c \in (0,1)$ such that for almost every environment  $\omega \in \Omega$ and for large $n \in \N$,  
\begin{equation} \label{eq_informal}
\mathbb{E}_{n, \omega}(\text{pm}(G_{n, \omega}(x))) \asymp n! c^n.
\end{equation}
Moreover, we have an explicit formula for the number $c$. This result implies that $\mathbb{E}_{n, \omega}(\text{pm}(G_{n, \omega}(x)))$ is a quenched variable.

The result in equation \eqref{eq_informal} should be understood in the sense that the mean number of perfect matchings for inhomogeneous random bipartite graphs in a random environment is asymptotically the same as the one of Erd\"os-R\'enyi graphs in which $p=c$. The number $c$ is the so-called \emph{scaling mean}  of a function related to the random graphs.  Scaling means were introduced, in more a general setting, in \cite{BIP_LLP} and are described in Sec.~\ref{sec:4} and Sec.~\ref{sec:5}.

\section{Inhomogeneous Random Bipartite Graphs in Random Environment} \label{sec:rbgre}
Consider the following generalization of the Erd\"os-R\'enyi graphs. Let  $W=\{w_1,  \dots , w_n\}$ and $M=\{m_1,  \dots, m_n\}$ be
two disjoint sets of  vertices. For every pair $1 \leq i,j \leq n$, let $a_{(ij)} \in [0,1]$ and consider the independent random variables 
 $X_{(ij)}$, with law
\begin{displaymath}
X_{(ij)} (x)= \left\{ \begin{array}{ll}
1 & \textrm{with probability } a_{(ij)} ;\\
0 & \textrm{with probability } 1-a_{(ij)}.
  \end{array} \right.
\end{displaymath}
Denote by  $G_n(x)$ the bipartite graph with  vertices $W, M$ and edges $E(x)$, where the edge $(w_i, m_j)$ belongs to $E(x)$ if and only if $X_{(ij)}(x)=1$. As it is clear from the definition all vertex of the graph do not play the same role. This contrasts with the (homogenous) Erd\"os-Renyi graphs (see \cite{BOLLOBAS2} for details).   In this note we consider inhomogoeneous random bipartite graphs in random environments, that is, the laws of $X_{(ij)}$ (and hence the numbers $a_{(ij)}$) are randomly chosen following an exterior environment law.  This approach to stochastic processes has  developed  since the ground breaking work by Solomon \cite{SOLOMON}  on Random Walks in Random Environment and subsequent work of a large community (see \cite{BOGACHEV} for a survey on the subject). 

The model we propose  is to consider the vertex sets $W, M$ as the environment and to consider that the number $a_{(ij)}$, which is the probability that the edge connecting $w_i$ with  $m_j$ occurs in the graph, is a random variable depending on   $w_i$ and $m_j$.
We remark that similar constructions has been studied in the setting of stochastic block model (see \cite{hll}). These have been used, for example, in machine learning in problems of community detection. We now describe precisely the graphs under consideration in this note.

The space of environments  is as follows. Fix $\alpha \in \N$ and  a stochastic vector $(p_1, p_2, \dots , p_{\alpha})$. Endow the set $\{1, \dots, \alpha\}$ with the probability  measure $P_{W}$ defined by $P_W(\{i\})=p_i$.  Denote by $\Omega_W$ the product space $\prod_{i=1}^{\infty}\{1, \dots , \alpha\}$ and by $\mu_{W}$ the corresponding product measure. For $\beta \in \N$ let $(\Omega_M, \mu_M)$ be the analogous probability measure space for the set $\{1, \dots, \beta\}$ and the stochastic vector $(q_1, q_2, \dots , q_{\beta})$. The \emph{space of environments} is $\Omega=\Omega_W\times \Omega_M$ with the measure $\mu_{\Omega}=\mu_W\times \mu_M$  and an \emph{environment} is an element $\omega \in \Omega$. Note that every environment defines two sequences
\begin{equation*}
W(\omega)=(w_1, w_2, \dots) \in \Omega_W \quad \textrm{and}\quad M(\omega)=(m_1, m_2, \dots) \in \Omega_M.
\end{equation*}

For each environment $\omega \in \Omega$ we now define the edge distribution  $X_{\omega, (ij)}$.  Let ${F}=[f_{sr}]$ be a $\alpha \times \beta$ matrix with entries $f_{sr}$ satisfying  $0\leq f_{sr}\leq 1$ and let 
$f:\{1,2 , \dots , \alpha\}  \times\{1,2 , \dots , \beta\} \to [0,1]$ be the function defined by $f(w, m)=f_{wm}$.  For each $\omega \in \Omega$ let
\begin{equation} \label{distribution2}
a_{(ij)}(\omega):=f\left(w_i(\omega), m_j(\omega)\right).
\end{equation}
Given an environment $\omega \in \Omega$ the corresponding \emph{edges distributions} are the random variables  $X_{\omega, (ij)}$
with laws
\begin{displaymath}
X_{\omega, (ij)} (x)= \left\{ \begin{array}{ll}
1 & \textrm{with probability } a_{(ij)}(\omega) ;\\
0 & \textrm{with probability } 1-a_{(ij)}(\omega).
  \end{array} \right.
\end{displaymath}

 Given an environment $\omega\in \Omega$, we  construct a sequence of random bipartite graphs $G_{n, \omega}$  considering the sets of vertices 
 \[
 W_{n, \omega}=(w_1(\omega), \dots, w_n(\omega))\quad \textrm{and}\quad M_{n, \omega}=(m_1(\omega), \dots, m_n(\omega)), 
 \]
  and edges distributions $X_{\omega, (ij)}$ given by the values of $a_{(ij)}(\omega)$ as in (\ref{distribution2}). We denote by $\mathbb{P}_{n, \omega}$  the law of the random graph $G_{n, \omega}$ and we call $F$ the \emph{edge distribution matrix}.

\begin{example} Given a choice of an environment $\omega\in \Omega$, the probability that the  bipartite graph $G_{n,\omega}(x)$ equals the complete bipartite graph $K_{n,n}$, using independence of the edge variables,  is
\[
\mathbb{P}_{n, \omega}\left(G_{n,\omega}(x)=K_{n,n}\right )=\prod_{1\leq i, j \leq n}\mathbb{P}_{n,\omega}(X_{\omega, (ij)}=1)=\prod_{1\leq i, j \leq n}a_{(ij)}(\omega).
\] 
\end{example} 

\section{Counting Perfect Matchings}\label{sec:4}
Recall that a perfect matching of a graph $G$ is a subset of edges containing every vertex exactly once. We denote by $\mathrm{pm}(G)$ the number of perfect matchings of $G$.  When the graph  $G$ is bipartite, and the corresponding partition of the vertices has the form $W=\{w_1, w_2, \dots , w_n\}$ and $M=\{m_1, m_2, \dots, m_n\}$, a perfect matching can be identified with a bijection between $W$ and $M$, and hence with a permutation $\sigma \in S_{n}$. From this, the total number of perfect matchings can be computed as
 \begin{equation}\label{eq.per}
 \mathrm{pm}(G)=\sum_{\sigma\in S_n}x_{1\sigma(1)}x_{2\sigma(2)}\cdots x_{n\sigma(n)}, 
 \end{equation}
where $x_{ij}$ are the entries of the adjacency matrix  $X_G$ of $G$, that is  $x_{ij}=1$ if  $(w_i, m_j)$ is an edge of $G$ and $x_{ij}=0$ otherwise. Of course, the right hand side of \eqref{eq.per} is the {\em permanent},  $\mathrm{per}(X_G)$, of the matrix $X_G$.  \\

In the framework of Section \ref{sec:rbgre}, we  estimate the number of perfect matchings for the sequence of inhomogeneous  random bipartite graphs $G_{n, \omega}$, for a given environment $\omega\in \Omega$. More precisely, we obtain estimates for the growth of the mean of
\begin{equation}\label{numpermat}
\text{pm}(G_{n, \omega}(x))=\per(X_{G_{n, \omega }(x)})=\sum_{\sigma \in S_n}X_{\omega, (1\sigma(1))}\cdots X_{\omega, (n\sigma(n))}.
\end{equation}
Denote by $\mathbb{E}_{n,\omega}$  the expected value with respect to the probability $\mathbb{P}_{n, \omega}$. Since the edges are independent and  $\mathbb{E}_{n, \omega}(X_{\omega, (ij)})=a_{ij}(\omega)$ we have
\begin{eqnarray*}\label{expnumpermat}
\mathbb{E}_{n, \omega}\left(\text{pm}(G_{n, \omega})\right)&=&\mathbb{E}_{n, \omega}\left(\sum_{\sigma \in S_n}X_{\omega, (1\sigma(1))}\cdots X_{\omega, (n\sigma(n))}\right) \\&=& \sum_{\sigma \in S_n}a_{ (1\sigma(1))}(\omega)\cdots a_{ (n\sigma(n))}(\omega)\\
&=&\per(A_{n} (\omega)), 
\end{eqnarray*} 
where the matrix $A_{n}( \omega)$ has entries $\left(A_{n}( \omega)\right)_{ij}=a_{(ij)}(\omega)$. The main result of this note describes the growth of this expected number for perfect matchings.

 The following number is a particular case of a quantity introduced by the authors in a more general setting in  \cite{BIP_LLP}.
\begin{definition} Let $F$ be an $\alpha \times \beta$  matrix with non-negative entries $(f_{rs})$. Let $\vec{p}=(p_1, \dots , p_{\alpha})$ and
$\vec{q}=(q_1, \dots , q_{\beta})$ be two stochastic vectors. The \emph{scaling mean} of $F$ with respect to $\vec{p}$ and $\vec{q}$ is defined by
\begin{equation*}
\text{sm}_{\vec{p}, \vec{q}}(F):= \inf_{(x_r) \in \R^{\alpha}_+, (y_s) \in \R^{\beta}_+} \left( \prod_{r=1}^{\alpha} x_r^{-p_r}\right)  \left( \prod_{s=1}^{\beta} y_s^{-q_s}\right)\left(\sum_{r=1}^{\alpha} \sum_{s=1}^{\beta} x_r f_{rs} y_s  p_rq_s\right).
\end{equation*}
\end{definition}

The scaling mean  
 is increasing with respect to the entries of the matrix and  lies between the minimum and the maximum of the entries  (see \cite{BIP_LLP} for details and more properties).
We stress that the scaling mean can be exponentially approximated using a simple iterative process (see Section \ref{sec:5}).
\smallskip

%
%

The main result in this note is the following,

\begin{theorem}[Main Theorem] \label{thm:main}
Let $(G_{n, \omega})_{n\geq 1}$ be a sequence of inhomogeneous random bipartite graphs in a random environment $\omega\in \Omega$. If the entries of the edge distribution matrix $F$ are strictly positive then  the following pointwise convergence holds
\begin{equation}\label{teorema}
\lim_{n\to \infty}\left(\frac{\mathbb{E}_{n, \omega}\left(\mathrm{pm}(G_{n, \omega})\right)}{n!}\right)^{1/n}=\mathrm{sm}_{\vec{p}, \vec{q}}(F),
\end{equation}
for $\mu_{W}\times \mu_{M}$-almost every environment $\omega\in \Omega$. 
\end{theorem}

\begin{remark}
As discussed in the introduction Theorem \ref{thm:main} shows that there exists a constant $c \in (0,1)$,   such that for almost every environment  $\omega \in \Omega$ and for $n \in \N$ sufficiently large
\begin{equation*} 
\mathbb{E}_{n, \omega}(\text{pm}(G_{n, \omega}(x))) \asymp n! c^n.
\end{equation*}
Namely, $c=\text{sm}_{\vec{p}, \vec{q}}(F)$.
This result should be compared with the corresponding one obtained by Erd\"os and R\'enyi for their class of random graphs, that is
\[\mathbb{E}(\text{pm}(G_n(x)))= n! p^n.\]
Thus, we have shown  that for large values of $n$ the growth of the number of perfect matchings for inhomogeneous  random graphs in a random environment behaves like the simpler model studied by Erd\"os and R\'enyi with $p=\text{sm}_{\vec{p}, \vec{q}}(F)$.
\end{remark}

\begin{remark}
Theorem \ref{thm:main} shows that the expected number of perfect matchings is a quenched variable.
\end{remark}

\begin{remark} Using the Stirling formula, the limit in (\ref{teorema}) can be stated as
\[
\lim_{n\to \infty}   \left( \frac{1}{n}\log\left(\mathbb{E}_{n, \omega}\left(\text{pm}(G_{n, \omega})\right)\right)-\log n \right)=\log \text{sm}_{\vec{p}, \vec{q}}(F)-1, 
\]
which gives a quenched result for the growth of the {\it perfect matching entropy} for the sequence of graphs $G_{\omega, n}$ (see \cite{ACFG}).
\end{remark}

\begin{remark}
Note that we assume a {\it uniform ellipticity} condition on the values of the probabilities $a_{(ij)}$ as in (\ref{distribution2}). A similar assumption appears in the setting of Random Walks in Random Environment (see \cite[p.355]{BOGACHEV}).
\end{remark}

We now present some concrete examples.
\begin{example}  \label{ex:1} Let $\alpha=\beta=2$ and $p_1=p_2=q_1=q_1= 1/2$. Therefore, the space of environments is the direct product of  two copies of the full shift on two symbols endowed with the $(1/2, 1/2)-$Bernoulli measure. 
The edge distribution matrix $F$ is a $2\times 2$ matrix with entries belonging to $(0,1)$. In \cite[Example 2.11]{BIP_LLP},  it was shown that
\begin{equation*} 
\text{sm}_{\vec p, \vec q} \begin{pmatrix}
f_{11} & f_{12} \\
f_{21} & f_{22}
\end{pmatrix} = \frac{\sqrt{f_{11}f_{22}} + \sqrt{f_{12}f_{21}}}{2} \, .
\end{equation*}
Therefore, Theorem \ref{thm:main} implies that
\begin{equation*}
\lim_{n\to \infty}\left(\frac{\mathbb{E}_{n, \omega}\left(\text{pm}(G_{n, \omega})\right)}{n!}\right)^{1/n}=\frac{\sqrt{f_{11}f_{22}} + \sqrt{f_{12}f_{21}}}{2},
\end{equation*}
for almost every environment $\omega\in \Omega$.
\end{example}

\begin{example}
 More generally let
$\alpha \in \N$ with $\alpha \geq 2$ and $\beta=2$. Consider the  stochastic vectors $\vec{p}=(p_1, p_2, \dots , p_{\alpha})$ and $\vec{q}=(q_1, q_2)$. The space  of environments is the direct product of a full shift on $\alpha$ symbols endowed with the  $\vec{p}$-Bernoulli measure with a full shift on two symbols endowed with the $\vec{q}$-Bernoulli measure. The edge distribution matrix $F$ is a $ \alpha \times 2$ matrix with entries $f_{r1}, f_{r2} \in (0, 1)$, where $r \in \{1, \dots , \alpha\}$.   Denote by $\chi \in \R^+$ the unique positive solution of the equation
\begin{equation*}
\sum_{r=1}^{\alpha} \frac{p_r f_{r1}}{f_{r1} + f_{r2} \chi} = q_1.
\end{equation*} 
Then
\begin{equation*} 
\text{sm}_{\vec{p},\vec{q}}(F)=
\text{sm}_{\vec p,\vec q} \begin{pmatrix}
f_{11} & f_{12} \\
\vdots & \vdots \\
f_{\alpha 1} & f_{\alpha 2}
\end{pmatrix} = q_1^{q_1} \left(\frac{q_2}{\chi}  \right)^{q_2}  \prod_{r=1}^{\alpha} \left(f_{r1} + f_{r2} \chi \right)^{p_r} \, .
\end{equation*}

Therefore, Theorem \ref{thm:main} implies that
\begin{equation*}
\lim_{n\to \infty}\left(\frac{\mathbb{E}_{n, \omega}\left(\text{pm}(G_{n, \omega})\right)}{n!}\right)^{1/n}= q_1^{q_1} \left(\frac{q_2}{\chi}  \right)^{q_2}  \prod_{r=1}^{\alpha} \left(f_{r1} + f_{r2} \chi \right)^{p_r},
\end{equation*}
for almost every environment $\omega\in \Omega$. The quantity in the right hand side first appeared in  work  by  Hal\'asz and Sz\'ekely in $1976$ \cite{HS},  in their study of  symmetric means.    In \cite[Theorem 5.1]{BIP_LLP} using a completely different approach we recover their result.
\end{example}

\section{Proof of the Theorem} The {\it shift map} $\sigma_{W}:\Omega_{W}\to \Omega_{W}$ is defined by 
\[
\sigma_{W}(w_1, w_2, w_3, \dots)=(w_2, w_3, \dots).
\]
The shift  map $\sigma_{W}$ is a $\mu_W$-preserving, that is,  $\mu_W(\Lambda)=\mu_W(\sigma_W^{-1}(\Lambda))$ for every measurable set $\Lambda\subset \Omega_W$,  and it is ergodic, that is,  if $\Lambda=\sigma_W^{-1}(\Lambda)$ then $\mu_W(\Lambda)$ equals $1$ or $0$. Analogously for $\sigma_M$ and $\mu_M$. We define a function $\Phi:\Omega_W\times \Omega_M\to \R$ by 
\[
\Phi(\vec w , \vec m)=f_{w_1 m_1}.
\]
Thus 
\[
\Phi(\sigma_W^{i-1}(\vec w), \sigma_M^{j-1}(\vec m))=f_{w_i  m_j}=a_{(ij)}(\omega).
\]
That is, the matrix $A_{n} (\omega)$ has entries $a_{(ij)}(\omega)=\Phi(\sigma_W^{i-1}(\vec w), \sigma_M^{j-1}(\vec m))$. We are in the  setting of the Law of Large Permanents  (\cite[Theorem 4.1]{BIP_LLP}).
\smallskip

\noindent{\bf Theorem (Law of Large Permanents).} 
Let $(X,\mu)$, $(Y,\nu)$ be Lebesgue probability spaces,
let $T \colon X \to X$ and $S \colon Y \to Y$ be ergodic measure preserving transformations,
and let $g\colon X \times Y \to \R$ be a positive measurable function essentially bounded away from zero and infinity.
Then for $\mu \times \nu$-almost every $(x, y) \in X \times Y$, the $n\times n$ matrix
$$
M_n(x, y)=\begin{pmatrix*}[l]
g(x,y)        & g(Tx,y)        & \cdots & g(T^{n-1}x,y)        \\
g(x,Sy)       & g(Tx,Sy)       & \cdots & g(T^{n-1}x,Sy)       \\
\qquad\vdots  & \qquad\vdots   &        & \qquad\vdots         \\
g(x,S^{n-1}y) & g(Tx,S^{n-1}y) & \cdots & g(T^{n-1}x,S^{n-1}y)
\end{pmatrix*}
$$ 

\[
\lim_{n\to \infty}\left(\frac{\per\left(M_n(x, y)\right)}{n!}\right)^{1/n}=\text{sm}_{\mu, \nu}(g)
\]
pointwise, where $\text{sm}_{\mu, \nu}(g)$ is the scaling mean of $g$ defined as
\[
\text{sm}_{\mu, \nu}(g)=\inf_{\phi, \psi} \frac{ \iint_{X\times Y}\phi(x)g(x,y)\psi(y)d\mu d\nu}{\exp\left(\int_{X} \log \phi(x)d\mu \right) \exp\left(\int_{Y} \log \psi(y)d\nu \right)},
\]
where the functions $\phi$ and $\psi$ are assumed to be measurable, positive and such that their logarithms are integrable.
\smallskip

Let $X=\Omega_W, Y=\Omega_M$, $T=\sigma_W, S=\sigma_M$, $g=\Phi$. As a consequence of an alternative characterization of the scaling mean (see \cite[Proposition 3.5]{BIP_LLP}) we have
\[
\text{sm}_{\mu_W, \mu_M}(\Phi)=\text{sm}_{\vec p, \vec q}(F).
\]
Since $f_{rs}>0$ we can apply the Law of Large Permanents to conclude the proof of the Main Theorem. $\hfill\blacksquare$

\begin{remark}
We have chosen to present our result in the simplest possible setting. That is, the environment space being products of full-shifts endowed with Bernoulli measures. In terms of stochastic block models we are considering only a finite number of communities. 
Using the general form of the Law of Large Permanent above  our results  can be extended  for inhomogeneous random graphs in  more general random environments.
\end{remark}

\section{A procedure to compute  the scaling mean}  \label{sec:5}
It is well known that  the computation of the permanent, and therefore of the number of perfect matchings,  is a very hard problem. Indeed, it was shown by Valiant \cite{v} that the evaluation of the permanent of $(0,1)$-matrices is an NP- hard problem. Thus, it is of interest to remark that the computation of the scaling mean, and therefore of the  expectation of the number of perfect matchings, 
can be performed with a simple iterative process that converges exponentially fast.

Denote by $\mathcal{B}^{\alpha} \subset \R^{\alpha}$ and by $\mathcal{B}^{\beta} \subset \R^{\beta}$  the positive cones. Define the following maps forming a (non-commutative) diagram:
$$
\begin{tikzcd}
\mathcal{B}^{\alpha} \arrow{r}{\mathtt{I}_1}  & \mathcal{B}^{\alpha} \arrow{d}{\mathtt{K}_2} \\
\mathcal{B}^{\beta} \arrow{u}{\mathtt{K}_1} & \mathcal{B}{\beta}\arrow{l}{\mathtt{I}_2}
\end{tikzcd}
$$
by the formulas:
\begin{eqnarray*}
(\mathtt{I}_1(\vec x))_i := \frac{1}{x_i}      \, ,  \quad
(\mathtt{I}_2(\vec y))_j  := \frac{1}{y_j}      \, ,&\\
(\mathtt{K}_2(\vec x))_j  := \sum_{i=1}^{\alpha} f_{ij}x_i p_i    \, , \quad
(\mathtt{K}_1(\vec y))_i  := \sum_{j=1}^{\beta} f_{ij}y_j q_j     .
\end{eqnarray*}
Let $\mathtt{T}: \mathcal{B}^{\alpha} \mapsto \mathcal{B}^{\alpha}$ be the map defined by
 $\mathtt{T} \coloneqq \mathtt{K}_1 \circ \mathtt{I}_2 \circ \mathtt{K}_2 \circ \mathtt{I}_1$. The map $\mathtt{T}$ is a contraction for a suitable Hilbert metric. Indeed, for $\vec x , \vec z \in \mathcal{B}^{\alpha}$ define the following (pseudo)-metric
\begin{equation*}
d(\vec x, \vec z):= \log \left( \frac{\max_i x_i/z_i}{\min_i x_i/z_i }		\right).
\end{equation*}
The following results were  proven in \cite[Lemma 3.4,  Lemma 3.3]{BIP_LLP}.
\begin{lemma}
For every $\vec{x}, \vec{z} \in \mathcal{B}^{\alpha}$  we have that
\begin{equation*}
d(T(\vec x), T(\vec z)) \leq \left( \tanh \frac{\delta}{4} \right)^2 d(\vec x, \vec z),
\end{equation*}
where
\[ \delta \leq 2 \log \left( \frac{\max_{i,j} f_{ij}}{\min_{ij } f_{ij}}		\right) < \infty.\]
\end{lemma}

 \begin{lemma} 
 The map $T$ has a unique (up to positive scaling) fixed point $\vec x_{\mathtt{T}} \in \mathcal{B}^{\alpha}$. Moreover, defining
 $\vec y_{\mathtt{T}} :=  \mathtt{K}_2 \circ \mathtt{I}_1(\vec x_{\mathtt{T}})$ one has that
 \begin{equation*}
 \text{sm}_{\vec{p}, \vec{q}}(F) = \prod^{\alpha}_{i=1} x_i^{p_i} \prod_{j=1}^{\beta} y_j^{q_j}.
 \end{equation*}
  \end{lemma}

Therefore, since the scaling mean can be directly computed from the fixed point of a contraction it possible to find good approximations of it using an iterative process that converges exponentially fast.

\begin{ack}
We appreciate interesting discussions with Manuel Cabezas about Random Walks in Random Environment.
\end{ack}


\bigskip

\bigskip

\begin{small}
	\noindent
	\begin{tabular}{lll}
		\textsc{Jairo Bochi} &
		\textsc{Godofredo Iommi} &
		\textsc{Mario Ponce} \\
		\email{\href{mailto:jairo.bochi@mat.puc.cl}{jairo.bochi@mat.puc.cl}} &
		\email{\href{mailto:giommi@mat.puc.cl}{giommi@mat.puc.cl}} &		
		\email{\href{mailto:mponcea@mat.puc.cl}{mponcea@mat.puc.cl}} \\
		\href{http://www.mat.uc.cl/~jairo.bochi}{www.mat.uc.cl/$\sim$jairo.bochi} &
		\href{http://www.mat.uc.cl/~giommi}{www.mat.uc.cl/$\sim$giommi} & 		
		\href{http://www.mat.uc.cl/~mponcea}{www.mat.uc.cl/$\sim$mponcea} 	
	\end{tabular}
	
	\bigskip
	
	\noindent
	\textsc{Facultad de Matem\'aticas, Pontificia Universidad Cat\'olica de Chile}
	
	\noindent
	\textsc{Avenida Vicu\~na Mackenna 4860, Santiago, Chile}
\end{small}

\end{document}